\documentclass[11pt]{amsart}
\usepackage[left=3.2cm,right=3.2cm,top=3.5cm,bottom=3.cm]{geometry}
\usepackage[T1]{fontenc}
\usepackage[utf8]{inputenc}
\usepackage[english]{babel}
\usepackage{amsmath}
\usepackage{amssymb}
\usepackage{amsthm}
\usepackage{braket}
\usepackage{amscd}
\usepackage{setspace}
\usepackage{booktabs}
\usepackage{lscape}
\usepackage{rotating}
\usepackage{tikz}
\usepackage{tikz-cd}
\usepackage{xcolor}
\usepackage{mathtools}
\usepackage{enumitem}
\usepackage{bigints}
\newtheorem{definition}{Definition}
\newtheorem{theorem}{Theorem}
\newtheorem*{theorem*}{Theorem}
\newtheorem{lemma}{Lemma}
\newtheorem{proposition}{Proposition}
\newtheorem{corollary}{Corollary}

\newcommand{\vare}{\varepsilon}
\newcommand{\Res}{\mathrm{Res}}
\newcommand{\Ext}{\mathrm{Ext}}

\newcommand{\Z}{\mathbb{Z}}
\newcommand{\C}{\mathbb{C}}

\newcommand{\R}{\mathbb{R}}
\newcommand{\pro}{\mathbb{P}}
\newcommand{\Po}{\mathbb{H}}

\newcommand{\SL}{\mathrm{SL}}
\newcommand{\GL}{\mathrm{GL}}
\newcommand{\Sym}{\mathrm{Sym}}

\def\={\;=\;}
\def\+{\;+\;}
\def\-{\;-\;}
\def\:={\;:=\;}

\title[Modular forms and deformation of punctures spheres]{Modular forms, deformation of punctured spheres, and extensions of symmetric tensor representations}
\author{Gabriele Bogo}
\address{Fachbereich Mathematik\\ Technische Universität Darmstadt\\ Schlossgartenstrasse 7, 64289 Darmstadt\\
Germany}
\email{bogo@mathematik.tu-darmstadt.de}

\begin{document}
\begin{abstract}
Let~$X=\Po/\Gamma$ be an~$n$-punctured sphere, $n>3$. We introduce and study~$n-3$ deformation operators on the space of modular forms~$M_*(\Gamma)$ based on the classical theory of uniformizing differential equations and accessory parameters. When restricting to modular functions, we recover a construction in Teichmüller theory related to the deformation of the complex structure of~$X$.
We describe the deformation operators in terms of derivations with respect to Eichler integrals of weight-four cusp forms, and in terms of vector-valued modular forms attached to extensions of symmetric tensor representations.
\end{abstract}

\maketitle

\section*{Introduction} 
Elliptic modular forms arise as solutions of linear differential equations related to the uniformization of hyperbolic Riemann surfaces. 
Consider for instance the weight one modular form~$f\in M_1(\Gamma_1(6))$
\begin{equation}
\label{eqn:fap}
f(q)\=1+3q+3q^2+3q^3+3q^4+\cdots\,,\quad q=e^{2\pi i\tau},\,\tau\in\Po
\end{equation}
that plays a role in Apery's proof of the irrationality of~$\zeta(2)$~\cite{beukers}. It is obtained as~$f(q)=y(t)\circ t(q)$ where~$y(t)$ is the holomorphic solution with~$y(0)=1$ of the Fuchsian differential equation
\begin{equation}
\label{eqn:apery}
\frac{d}{dt}\Bigl(t(t-1)(t-1/9)\frac{d}{dt}\Bigr)y(t)+(t-\rho)y(t)\=0\,
\end{equation}
with~$\rho=1/3$, and~$t(q)$ is a Hauptmodul for~$\Gamma_1(6)$.
The differential equation~\eqref{eqn:apery} is associated to the uniformization of the four-punctured sphere~$\pro^1\smallsetminus\{\infty,1,0,1/9\}$; the parameter~$\rho$ is called~\emph{accessory parameter}. 

For every choice of the accessory parameter~$\rho$, let~$\{y_\rho(t),\hat{y}_\rho(t)\}$ be a Frobenius basis of solutions of~\eqref{eqn:apery}. If~$Q_\rho:=\exp(\hat{y}_\rho(t)/y_\rho(t))$, one can construct a~$Q$-expansion of the form
\begin{equation}
\label{eqn:esdef}
f_\rho(Q)\=y_\rho(t)\circ t_\rho(Q),\quad t_\rho(Q):=Q_\rho(t)^{-1}\,
\end{equation}
where~$t_\rho(Q)$ is the compositional inverse of~$Q_\rho(t)$. A consequence of uniformization theory is that when~$\rho=1/3$ the function~$f_\rho(Q)$ is precisely~$f(q)$ in~\eqref{eqn:fap}.

We make use of the dependence of~$f_q(Q)$ on the accessory parameter~$\rho$ to define a \lq\lq{deformation}\rq\rq\;(the choice of the name will be explained later) of the Fourier expansion of~$f(q)$ by
\begin{equation}
\label{eqn:defes}
\partial_0{f}(q)\:=\frac{\partial f_\rho(Q)}{\partial\rho}\Bigl|_{\rho=1/3}\=9q+\frac{153}{2}q^2+105q^3+\frac{543}{4}q^4+\frac{36057}{200}q^5-\frac{17607}{200}q^6+\cdots\,.
\end{equation}
By looking at the above coefficients one realizes that~$\partial_0{f}$ can be written as
\begin{align*}
\partial_0{f}(q)&\=\Bigl(1+3q+3q^2+3q^3+\cdots\Bigr)\Bigl(9q-9\frac{q}{2}^2-3q^3+\frac{9}{4}q^4+\cdots\Bigr)\\
&\+2\Bigl(9q-\frac{9}{4}q^2-q^3+\frac{9}{16}q^4+\cdots\Bigr)\Bigl(3q+6q^2+9q^3+12q^4+\cdots\Bigr)\\
&\=f(q)\widetilde{h}'(q)+2\widetilde{h}(q)f'(q)\,
\end{align*}
where~$h(q)=9q-18q^2-27q^3+\cdots\in S_4(\Gamma_1(6))$,~$\widetilde{h}(q)=9q-9/2q^2-3q^3+\cdots $ is the Eichler integral of~$h$, and~$'=q\tfrac{d}{dq}$.

\smallskip

More generally, let~$X=\Po/\Gamma$ be an~$n$-punctured sphere,~$n\ge3$. The classical theory of uniformization attaches to~$X$ a family of second-order Fuchsian differential equations depending on~$n-3$ \emph{accessory parameters}~$\rho_0,\dots,\rho_{n-4}$. As in~\eqref{eqn:esdef}, one can construct a~$Q$-expansion~$f_\rho(Q)$ for any choice of the accessory parameters~$\rho=(\rho_0,\dots,\rho_{n-4})$. As follows from the uniformization theorem, there exists a unique value~$\rho_F$ of these parameters, called the \emph{Fuchsian value}, that makes~$f_{\rho_F}(q)$ the Fourier expansion at a cusp of a modular form on the uniformizing Fuchsian group~$\Gamma$.
Analogously to~\eqref{eqn:defes}, for every $i=0,\dots,n-4$ one can consider the derivative~$\partial{f_\rho(Q)}/\partial\rho_i$ and then specialize~$\rho$ to the Fuchsian value~$\rho_F$ in order to define a new holomorphic function on~$\Po$. This operation extends to a differential operator~$\partial_i$ on the space of modular forms~$M_*(\Gamma)$.

\smallskip

In this paper we study the deformation operators~$\partial_i$ and characterize them in three ways. We summarize our results in the following theorem.
\begin{theorem*}
Let~$\Po/\Gamma$ be isomorphic to an~$n$-punctured sphere~$n\ge4$, and let~$g\in{M}_k(\Gamma)$. 
\begin{enumerate}
\item There exist a basis~$\{h_0,\dots,h_{n-4}\}$ of~$S_4(\Gamma)$ such that
\[
{\partial_i{g}}\=kg\widetilde{h}_i'\+2\widetilde{h}_ig'\=[g,\widetilde{h}_i]_1\,,\quad i=0,\dots,n-4\,,
\]
where~$\widetilde{h}_i$ is the Eichler integral of the cusp form~$h_i$ and~$[\;,\,]_1$ is the first Rankin-Cohen bracket.
\item Let~$t$ be a Hauptmodul for~$\Gamma$ and let~$\nu_i$ be the harmonic Beltrami differential dual to the cusp form~$h_i$. Let~$\vare>0$ be such that~$\|\vare\nu_i\|_\infty<1$ and let~$\Gamma^{\vare\nu_i}$ be the  group obtained by conjugating~$\Gamma$ by a quasiconformal solution of the Beltrami equation associated to~$\vare\nu_i$.  If~$t^{\vare\nu_i}$ is a Hauptmodul for~$\Gamma^{\vare\nu_i}$, then
\[
\partial_i t\=4\,\frac{\partial{t^{\vare\nu_i}}}{\partial\bar{\vare}}\Bigr|_{\vare=0}.
\]
\item The~$i$-th deformation operator~$\partial_i$ induces a map
\[
\overset{\rightarrow}{\partial}_i\colon M_k(\Gamma)\to M_k(\Gamma,v^{h_i}_{0,2})\,
\]
from the space of weight~$k$ modular forms to the space of weight~$k$ vector-valued modular forms with respect to an extension~$v^{h_i}_{0,2}$ of symmetric tensor representations~$v_0$ and~$v_2$ of dimension~$1$ and~$3$ respectively.
\end{enumerate}
\end{theorem*}

In Section~\ref{sec:pre} we collect some basic facts on classical uniformization and modular forms, symmetric tensor representations, their extensions, and vector-valued modular forms. 
In Section~\ref{sec:def} we define the deformation operators. Part 1 of the theorem is proven in Sections~\ref{sec:defde} and~\ref{sec:defmf}; the proof is based on the theory of differential equations. 
Section~\ref{sec:teich} contains the proof of the second statement and some background in Teichmüller theory. The result on vector-valued modular forms, which follows from the proof of the first point, can be found in Section~\ref{sec:vvmf}. The last section~\ref{sec:final} contains some final remarks and open questions.

\section{Uniformization, differential equations, symmetric tensor representations}
\label{sec:pre}
In this section we recall some basic facts and fix notation.
\subsection{Uniformization and accessory parameters}
\label{sec:premod}
To an~$n$-punctured sphere~$X=\pro^1\smallsetminus\{a_1,\dots,a_{n-1}=0,a_n=\infty\},\,n\ge 3,$ one can attach a family of second-order Fuchsian differential equations
\begin{equation}
\label{eqn:de}
L_X\:=\frac{d}{dt}\Bigl(P(t)\frac{d}{dt}\Bigr)\+\sum_{i=0}^{n-3}{\rho_it^i}\,,\quad P(t)\:=\prod_{j=1}^{n-i}(t-a_j)\,
\end{equation}
depending on~$n-3$ parameters~$\rho_0,\dots,\rho_{n-4}$ called~\emph{accessory parameters}. The value of the parameter~$\rho_{n-3}=(n/2-1)^2$ is fixed to make the singular point at~$\infty$ regular singular and for every choice of~$\rho_1,\dots,\rho_{n-3}$ all finite singular points are regular with local exponents~$(0,0)$.

This family of differential equations is classically related to the Fuchsian uniformization of~$X$~\cite{sgervais}. From the uniformization theorem it follows that there is a unique choice of the parameters~$\rho=(\rho_0,\dots,\rho_{n-4})$ such that the ratio of linearly independent solutions of~\eqref{eqn:de} lifts to a biholomorphic map between the universal covering of~$X$ and the upper half-plane~$\Po$. This map gives a universal covering of~$X$. We call this special choice of parameters the~\emph{Fuchsian value} and denote it~$\rho_F$. 
In this special case, the monodromy group~$\Gamma$ of~\eqref{eqn:de} is the Deck group of the universal covering~$\Po\to X$. This implies that~$\Gamma\subset\SL_2(\R)$ is a Fuchsian group, i.e., discrete and cofinite. A consequence is that a holomorphic solution of~\eqref{eqn:de} in the case~$\rho=\rho_F$ lifts to a holomorphic function on~$\Po$ that is a~($k$-th root of a weight~$k$) modular form on~$\Gamma$. More details and examples can be found in Chapter 5 of Zagier's exposition in~\cite{123}.

\smallskip

For every choice of~$\rho$ in~\eqref{eqn:de} one can construct some power series from a Frobenius basis of solutions~$\{y_\rho(t),\hat{y}_\rho(t)\}$ at~$t=0$:
\begin{align}
\label{eqn:fncts}
Q_\rho(t)&\:=\exp(\hat{y}_\rho(t)/y_\rho(t))\=\sum_{s=1}^\infty{Q_s(\rho)t^s}\,,\\
t_\rho(Q)&\:=Q_\rho(t)^{-1}\=\sum_{s=1}^\infty{t_s(\rho)Q^s}\,,\\ 
f_\rho(Q)&\:=y_\rho(t_\rho(Q))\=\sum_{s=0}^\infty{f_s(\rho)Q^s}\,,
\end{align}
where~$t_\rho(Q)$ is the compositional inverse of~$Q_\rho(t)$. In all the above power series the coefficients are polynomials in the accessory parameters. As recalled in the previous paragraph, if~$\tau\in\Po$,  when~$\rho=\rho_F$ we have that~$\hat{y}_{\rho_F}(t)/y_{\rho_F}(t)=2\pi i \tau+\lambda$ for some~$\lambda\in\C$. Consequently
\[
Q_{\rho_F}(t)\=e^{2\pi i\tau+\lambda}\=c\cdot{q},\quad q=e^{2\pi i\tau}\,,\tau\in\Po\,,
\]
and
\begin{equation}
\label{eqn:tf}
t(\tau):=t_{\rho_F}(cq),\quad f(\tau):=f_{\rho_F}(cq)
\end{equation}
are holomorphic functions on~$\Po$. More precisely,~$t(\tau)$ is a Hauptmodul for~$\Gamma$ and~$f(\tau)$ is a root of a modular form. As shown in the appendix of~\cite{bogo}, it turns out that~$f$ is the square root of a weight two modular form with all its zeros concentrated in the cusp~$x$ where the Hauptmodul~$t$ has its unique pole.


In the following we will assume that~$f$ is itself a modular form, i.e., that~$f\in M_1(\Gamma)$; this is not always the case, for instance  if~$-I\in\Gamma$, but this assumption makes the exposition smoother and we do not lose much in terms of generality. The reader can check all the statements in the next sections can easily be adapted to the more general case~$f^2\in M_2(\Gamma)$.

\subsection{Symmetric tensor representations and extensions}
\label{sec:ext}
Closely related to the differential equations solved by modular form are certain representations of the Fuchsian group~$\Gamma$.
Let~$r\ge 0$ and~$V_r:=\Sym^{r}(\C^2)$ and let~$v_r\colon\SL_2(\R)\to \mathrm{GL}(V_r)$ be the representation defined by
\[
v_r(\gamma)\begin{pmatrix} z_1 \\ z_2\end{pmatrix}^r\=\biggl(\gamma\begin{pmatrix} z_1 \\z_2\end{pmatrix}\biggr)^r\,
\]
for every~$\gamma\in\SL_2(\R),$ where
\[
\begin{pmatrix} z_1 \\ z_2\end{pmatrix}^r:={}^t\bigl(z_1^r,z_1^{r-1}z_2,\dots,z_1z_2^{r-1},z_2^r\bigr)\in\C^{r+1}\,.
\]
We denote by~$M_k(\Gamma,v_r),\,S_k(\Gamma,v_r)$ the space of holomorphic vector-valued modular forms (VVMFs for short) and vector-valued cusp forms with respect to the representation~$v_r$. These VVMFs are strictly related to quasimodular forms; as Kuga-Shimura~\cite{KS} and Choie-Lee~\cite{CL} showed, there is an isomorphism
\begin{equation}
\label{eqn:cl}
M_k(\Gamma,v_r)\;\simeq\;\oplus_{s=0}^r{M_{k+r-2s}(\Gamma)}\;\simeq\;\widetilde{M}_{r+k}(\Gamma)^{(\le r)}\,,\quad\text{for every }k,r\ge 0\,,
\end{equation}
where~$\widetilde{M}_{r+k}(\Gamma)^{(\le r)}$ denotes the space of quasimodular forms of weight~$r+k$ and depth~$\le r$.
The reader may consult~\cite{CL} for a proof of this result involving Rankin-Cohen brackets, and as a reference for symmetric tensor representations and modular forms.

\smallskip

We will be interested in extensions of symmetric tensor representations, i.e., elements of~$\Ext^1_\Gamma(v_l,v_r)$ for some~$l,r\ge0$. We can describe the space~$\Ext^1_\Gamma(v_l,v_r)$ in terms of quasimodular forms and VVMFs as follows.
From the identity~$\Ext^1_\Gamma(v_l,v_r)=H^1(\Gamma,v_l\otimes v_r)$, together with the Clebsch-Gordan decomposition of~$v_l\otimes v_r$ and the classical Eichler-Shimura isomorphism one can prove that, if~$r\le l$,
\begin{equation}
\label{eqn:ext}
0\longrightarrow{\bigoplus_{s=0}^{r}M_{l+r+2-2s}(\Gamma)}\longrightarrow{\Ext^1_\Gamma(v_l,v_r)}\longrightarrow{\bigoplus_{s=0}^{r}S_{l+r+2-2s}(\Gamma)}\longrightarrow0
\end{equation}
is a short exact sequence. Together with~\eqref{eqn:cl} this implies that
\[
0\longrightarrow{M_{r+2}(\Gamma)}\longrightarrow{\Ext^1_\Gamma(v_l,v_r)}\longrightarrow{S_{r+2}(\Gamma)}\longrightarrow0\,
\] 
is exact. 

If~$v_{l,r}$ is an extension of~$v_r$ by~$v_l$ we denote by~$M_k(\Gamma,v_{l,r})$ the space of holomorphic modular forms with respect to the representation~$v_{l,r}$. Examples of these VVMFs have been studied, among other things, in a recent paper by Mertens and Raum~\cite{MR}.

\section{Deformation operators}
\label{sec:def}
Let~$\Gamma\subset\SL_2(\R)$ be a torsion-free genus zero Fuchsian group, $t\colon\Po/\Gamma\to X$ a Hauptmodul and let~$f\in M_1(\Gamma)$ be the modular solution of the differential equation \eqref{eqn:de} in the case~$\rho=\rho_F$. Recall from the last section that~$f\=f_\rho(q)|_{\rho=\rho_F},\, t=t_\rho(q)|_{\rho=\rho_F}$ where
\[
f_\rho(q)=\sum_{n=0}^\infty{f_n(\rho)q^n},\quad t_\rho(q)=\sum_{n=0}^\infty{t_n(\rho)q^n},\quad\rho=(\rho_0,\dots,\rho_{n-4})
\]
and~$\rho_F$ is the Fuchsian parameter. For~$i=0,\dots,n-4$ we consider the partial derivatives
\[
\frac{\partial{f_\rho(q)}}{\partial\rho_i}\:=\sum_{n=0}^\infty{\frac{\partial f_n(\rho)}{\partial\rho_i}q^n},\quad \frac{\partial{t_\rho(q)}}{\partial\rho_i}\:=\sum_{n=0}^\infty{\frac{\partial t_n(\rho)}{\partial\rho_i}q^n}\,.
\]
Recall from the discussion in Section~\ref{sec:premod} that every~$g\in M_k(\Gamma)$ is of the form~$g=f^kR(t)$ where~$R\in\C[t]$.
\begin{definition}
Let~$k\ge 0$ and~$g\in M_k(\Gamma),\, g=f^kR(t)$. For~$i=0,\dots,n-4$, define the~\emph{$i$-th deformation}~$\partial_i{g}$ of~$g$ by
\[
\partial_i{g}(q)\:=\frac{\partial f_\rho(q) R\bigl(t_\rho(q)\bigr)}{\partial\rho_i}\biggl|_{\rho=\rho_F}\,.
\]
\end{definition}
For every~$i=0,\dots,n-4,$, the~$i-$th deformation operator~$\partial_i$ defines a derivation on the space of modular forms. It can be described in terms of the first Rankin-Cohen bracket and Eichler integrals of weight four cusp forms as follows.
\begin{theorem}
\label{thm:main}
Let~$\Po/\Gamma$ be isomorphic to an~$n$-punctured sphere, and let~$g\in{M}_k(\Gamma)$. There exist a basis~$\{h_0,\dots,h_{n-4}\}$ of~$S_4(\Gamma)$ such that
\[
{\partial_i{g}}\=kg\widetilde{h}_i'\+2\widetilde{h}_ig'\=[g,\widetilde{h}_i]_1\,,\quad i=0,\dots,n-4\,,
\]
where~$\widetilde{h}_i$ denotes the Eichler integral of the cusp form~$h_i$.
\end{theorem}
The proof of the theorem is given in the next two subsections. We study first the effect of the partial derivation~$\partial/\partial\rho_i$ on the solution~$y_\rho(t)$ of the differential equation~\eqref{eqn:de}; it turns out that~$\partial{y_\rho(t)}/\partial\rho_i$ satisfies an higher-order Fuchsian differential equation whose shape depends on~$i$. The solutions of these higher-order ODEs have an integral representation in terms of the function~$y_\rho(t)$. Second, we show that these solutions give rise to Eichler integrals of cusp forms (and higher iterated integrals of modular objects) as functions on the universal covering~$\Po$.

\subsection{Deformation of differential equations}
\label{sec:defde}
\begin{proposition}
\label{thm:defde}
Let~$\{y_\rho(t),\hat{y}_\rho(t)\}$ be a fundamental system of solution of~\eqref{eqn:de}. Then for every~$i=0,\dots,n-4$ and $r\ge0,$ the functions
\[
\frac{\partial^ry_\rho(t)}{\partial\rho_i^r}\Bigl|_{\rho=\rho_F},\,\frac{\partial^r\hat{y}_\rho(t)}{\partial\rho_i^r}\Bigl|_{\rho=\rho_F}
\]
are solutions of an order~$2(r+1)$ Fuchsian differential equation with the same singular points as~\eqref{eqn:de}. The associated monodromy representation~$\colon\pi_1(X,x_0)\to\GL_{2(r+1)}(\C)$ is an iterated extension of representations~$\pi_1(X,x_0)\to\Gamma\subset\SL_2(\R)$.
\end{proposition}

The first claim of the theorem is independent on the choice of the accessory parameters. In fact, a more general result holds.
\begin{proposition}
\label{thm:comp}
Let~$\{y_\rho,\hat{y}_\rho\}$ be a fundamental system of solutions of~\eqref{eqn:de}. Then for every~$i=0,\dots,n-4$ and~$r\ge0$ the functions~$\partial^ry_\rho(t)/\partial\rho_i^r$ and~$\partial^r\hat{y}_\rho(t)/\partial\rho_i^r$ satisfy an order~$2(r+1)$ Fuchsian differential equation with the same singular points as~\eqref{eqn:de}.
\end{proposition}
\begin{proof}
We work with a Frobenius basis of solutions~$\{y_\rho(t),\hat{y}_\rho(t)=\log(t)y_\rho(t)+\tilde{y}_\rho(t)\}$ near the regular singular point~$t=0$.
Consider the deformed Fuchsian operator
\begin{equation}
\label{eqn:def}
L_{X,\epsilon}\:=\frac{d}{dt}\Bigl(P(t)\frac{d}{dt}\Bigr)\+\sum_{i=0}^{n-3}{(\rho_i+\epsilon_i)t^i}
\end{equation}
where~$\rho_{n-3}=(1-n/2)^2$ as before and~$\epsilon_{n-3}=0$. We do not deform~$\rho_{n-3}$ in order to preserve the regular singularity at~$\infty$.

\smallskip

Let~$m:=n-4$ and let
\[
A_s(\rho,\epsilon):=\sum_{j\in\Z_{\ge0}^m}{A_s^{(j)}(\rho)\epsilon^j},\qquad \epsilon^j:=\epsilon_0^{j_0}\cdots\epsilon_m^{j_m},
\] 
be the solution of the homogeneous recursion associated to the differential operator~$L_{X,\epsilon}$ with initial values~$(A_{3-n},\dots, A_{-1},A_0)=(0,\dots,0,1)$. Finally, consider the power series
\[
\Phi_j(\rho,t)\:=\sum_{s=0}^\infty{A^{(j)}_s(\rho)t^s},\qquad \Phi(\rho,\epsilon,t)\:=\sum_{j\in\Z^m_{\ge0}}{\Phi_j(\rho,t)\epsilon^j}\,.
\]
We have the following identity
\begin{equation}
\label{eqn:e=p}
\biggl(\frac{\partial}{\partial\epsilon_i}-\frac{\partial}{\partial\rho_i}\biggr)\Phi(\rho,\epsilon,t)\=0\,,
\end{equation}
which can be proved from the definition of the deformed operator~$L_{X,\epsilon}$. In fact, the linear recursion associated to~\eqref{eqn:def} is of the form
\begin{equation}
\label{eqn:rec}
\hat{P}_{n-2}(s)A_{s+1}(\rho,\epsilon)\=\sum_{i=0}^{n-3}{P_i(s,\rho,\epsilon)A_{s-i}(\rho,\epsilon)},\quad P_i(s,\rho,\epsilon):=\rho_i+\epsilon_i+\hat{P}_i(s)
\end{equation}
where~$\hat{P}_i(s),\,i=0,\dots,n-2$ is a quadratic polynomial in~$s$ that does not depend on~$\rho,\epsilon$ and~$\hat{P}_{n-2}(s)=(s+1)^2\prod_{j=1}^{n-2}{(-a_j)}$. By induction on~$s$ we see that
\[
\hat{P}_{n-2}(s)\frac{\partial A_{s+1}(\rho,\epsilon)}{\partial\rho_{i_0}}\=A_{s-i_0}(\rho,\epsilon)+\sum_{i=0}^{n-3}{P_i(s,\rho,\epsilon)\frac{\partial A_{s-i}(\rho,\epsilon)}{\partial \rho_{i_0}}}\=\hat{P}_{n-2}(s)\frac{\partial A_{s+1}(\rho,\epsilon)}{\partial\epsilon_{i_0}}\,,
\]
which immediately implies~\eqref{eqn:e=p}. In particular, we get 
\[
\frac{\partial\Phi_{(j_0,\dots,j_i,\dots,j_m)}(\rho,t)}{\partial\rho_i}\=(j_i+1)\Phi_{(j_0,\dots,j_i+1,\dots,j_m)}(\rho,t)\,.
\]
By applying the above identity~$r$ times to the function~$\Phi_{(0,\dots,0)}(\rho,t)=y_\rho(t)$ we find
\[
\frac{\partial^ry_\rho(t)}{\partial\rho_i^r}\=\frac{\partial^r\Phi_0(\rho,t)}{\partial\rho_i^r}\=(r+1)!\Phi_{(0,\dots0,r,0,\dots,0)}(\rho,t)\,.
\] 
This shows that~$\partial/\partial\rho_i$ acts on~$\Phi_j(\rho,t)$ by raising the~$i$-th index.

On the other hand, the Fuchsian operator~$L_X$ acts by lowering the indices of~$\Phi_j(\rho,t).$ In fact, from the definition of~$L_{X,\epsilon}$ we see that
\begin{equation}
\label{eqn:lphi}
L_X\Phi(\rho,\epsilon,t)\=-\Phi(\rho,\epsilon,t)\sum_{i=0}^{m}{\epsilon_it^i}\,,
\end{equation}
which implies
\begin{equation}
\label{eqn:lact}
-L_X\Phi_{(j_0,\dots,j_k,\dots,j_m)}(\rho,t)\=\Phi_{(j_0-1,\dots,j_k,\dots,j_m)}(\rho,t)\+\cdots+t^m\Phi_{(j_0,\dots,j_k,\dots,j_m-1)}(\rho,t)\,.
\end{equation}
The idea is to contrast the action of~$\partial/\partial\rho_i$ (raising indices) by applying the differential operator~$L_X$ (lowering indices) but, as we see from \eqref{eqn:lact}, some powers of~$t$ appears. To deal with this, we introduce new second-order differential operators~$L_{X,r},\,r\ge1$ by
\begin{equation}
\label{eqn:lk}
L_{X,r}\bigl(t^rY\bigr)\=t^{r+1}L_X\bigl(Y\bigr),\quad r\ge1\,.
\end{equation}
It is easy to see that~$L_{X,r}$ is a Fuchsian operator with the same singular points as~$L$ for every~$r\ge 1$ (but not with the same local exponents).
It follows immediately from~\eqref{eqn:lk} and~\eqref{eqn:lphi} that~$L_{X,r}\bigl(t^r\Phi(\rho,\epsilon,t)\bigr)=-\Phi(\rho,\epsilon,t)\sum_{i=0}^{m}{\epsilon_it^{i+r+1}}$
and in particular
\begin{equation}
\label{eqn:lkact}
\frac{-L_{X,r}\bigl(t^r\Phi_{(j_0,\dots,j_r,\dots,j_m)}(\rho,t)\bigr)}{t^{r+1}}\=\Phi_{(j_0-1,\dots,j_r,\dots,j_m)}(\rho,t)\+\cdots+t^{m}\Phi_{(j_0,\dots,j_r,\dots,j_m-1)}(\rho,t)\,.
\end{equation}
If~$\Phi_{r_i}(\rho,t)$ denotes the function~$\Phi_{(0,\dots,0,r,0,\dots,0)}$ where the unique non-zero index~$r$ is at the~$i$-th place, 
from~\eqref{eqn:lact} and~\eqref{eqn:lkact} we conclude that
\[
\begin{cases}
\bigl(L_X^{r+1}\bigr)\Phi_{r_i}(\rho,t)&\=0\quad i=0\,,\\
\bigl(L_{X,r(i+1)-1}\circ\cdots\circ L_{X,2i+1}\circ L_{X,i}\circ L_X\bigr)\Phi_{r_i}(\rho,t)&\=0\quad i=1,\dots,m\,,
\end{cases}
\]
from which the statement for~$\partial^ry_\rho(t)/\partial\rho_i^r=\Phi_{r_i}(\rho,t)$ follows.

\smallskip

The same argument, with the same differential operators, works for~$\partial^r\hat{y}_\rho(t)/\partial\rho_i^r$ by deforming the recursion giving the coefficients of the holomorphic part~$\tilde{y}_\rho(t)$ of~$\hat{y}_\rho(t)$.
\end{proof}

In order to discuss the monodromy, we need the following lemma.
\begin{lemma}
\label{thm:ext}
Let~$M_1,M_2$ be linear differential operators with the same set of singular points in~$\pro^1$ and let~$M:=M_2\circ M_1\,.$
The monodromy representation associated to a fundamental system of solutions of~$MY=0$ is an extension of monodromy representations associated to solutions of~$M_1Y=0$ and~$M_2Y=0$ respectively. 
\end{lemma}
\begin{proof}
Let~$D:=\pro^1\smallsetminus\{\text{singular points of }M_1\}$. 
Without loss of generality, we can choose a fundamental system of solution~$Y_M$ of~$MY=0$ of the form
\[
Y_M=\bigl(u_0,\dots,u_k,v_1,\dots,v_l\bigr)
\]
where~$Y_{M_1}=(u_0,\dots,u_k)$ and~$Y_{M_2}=(M_1v_1,\dots,M_1v_l)$ are fundamental system of solutions of~$M_1Y=0$ and~$M_2Y=0$ respectively.
If~$\rho_*\colon\pi_1(D,d_0)\to\GL_{**}(\C)$ denotes the monodromy representation associated to~$Y_*,\,*=M,M_1,M_2$, it is easy to verify that the sequence of~$\pi_1(D,d_0)$-modules
\[
\begin{tikzcd}
0\arrow[r] & (S_{M_1},\rho_1)\arrow[r,"\phi"] &(S_M,\rho)\arrow[r,"\psi"]&(S_{M_2},\rho_2)\arrow[r] &0\,,
\end{tikzcd}
\]
where~$S_*=\mathrm{Span}(Y_*)$, is exact.
\end{proof}

\begin{proof}[proof of Proposition 1]
We need only to prove the statement about the monodromy representation.
When~$\rho=\rho_F$ a basis of solutions of $L_XY=0$ is given locally on the universal covering~$\Po$ by~$\{\sqrt{f},\tau\sqrt{f}\}$, where~$f\in M_2(\Gamma)$ is a weight two modular form and~$\Gamma\subset\SL_2(\R)$ is the Deck group of the covering~$\Po\to X$.
It follows from the definition~\eqref{eqn:lk} that for every~$r\ge1$ a basis of solutions of~$L_{X,r}Y=0$ is given locally on~$\Po$ by~$\{t^r\sqrt{f},\tau t^r\sqrt{f}\}$, i.e., that~$L_{X,r}Y=0$ has (meromorphic) modular solutions. 

By Proposition~\ref{thm:comp}, $\partial^ry_\rho(t)/\partial\rho_i^r,\,i=0,\dots,n-4$, is annihilated by a composition of the differential operators~$L_X$ and $L_{X,j},\,j=1,\dots, r(i+1)-1$. It follows from Lemma~\ref{thm:ext} that the monodromy of the associated differential equation is an (iterated) extension of monodromy representations of~$L_X,\,L_{X,j}$. It is a general fact that the monodromy group of a second-order differential equation with modular solutions in~$M_*(\Gamma)$ is conjugated to~$\Gamma$ (see Chapter 5 of Zagier's exposition in~\cite{123}).
\end{proof}

\begin{corollary}
\label{thm:integralform}
The following recursive formula holds for~$r\ge0,\,i=0,\dots,n-4$:
\[
\frac{\partial^{r+1}y_\rho(t)}{\partial\rho_i^{r+1}}\=y_\rho(t)\!\!\bigintsss_0^t{\!\!\frac{1}{y_\rho^2(t)P(t)}{\bigintssss_{0}^{t_1}\!\!\!y_\rho(t_2)\frac{\partial^{r}y_\rho(t_2)}{\partial\rho_i^{r}}\;dt_2}\;dt_1}\,. \\
\]
\end{corollary}
\begin{proof}
For every holomorphic function~$u$ one has
\begin{equation}
\label{eqn:intrep}
L_X\Biggl(y_\rho(t)\!\!\bigintsss_0^t{\!\!\frac{1}{y_\rho^2(t)P(t)}{\bigintssss_{0}^{t_1}\!\!\!y_\rho(t_2)u(t_2)\;dt_2}\;dt_1}\Biggr)\=u\,,
\end{equation}
as follows from a straightforward computation. 
From the definition~\eqref{eqn:lk} it follows moreover that, if~$i>0$,
\begin{equation}
\label{eqn:rlk}
L_{X,r(i+1)-1}\circ\cdots\circ L_{X,2i+1}\circ L_{X,i}\circ L_X(u)\=t^{r(i+1)}L_X\circ\underbrace{\frac{L_X}{t^i}\circ\frac{L_X}{t^i}\circ\cdots\circ\frac{L_X}{t^i}}_{r-1}(u)\,.
\end{equation}
Now let~$v_i(t)$ be such that~$L_{X,r(i+1)-1}\circ\cdots\circ L_{X,2i+1}\circ L_{X,i}\circ L_X(v_i)=0$. Then the function
\begin{equation}
\label{eqn:wi}
w_i(t)\:=y_\rho(t)\!\!\bigintsss_0^t{\!\!\frac{1}{y_\rho^2(t)P(t)}{\bigintssss_{0}^{t_1}\!\!\!y_\rho(t_2)t^iv_i(t_2)\;dt_2}\;dt_1}
\end{equation}
is such that~$L_{X,(r+1)(i+1)-1}\circ\cdots\circ L_{X,i}\circ L_X(w_i)=0$ if~$i\ge 1$. We have in fact
\begin{align*}
L_{X,(r+1)(i+1)-1}\circ\cdots\circ L_{X,i}\circ L_X(w_i)&\= t^{(r+1)(i+1)}L_X\circ\underbrace{\frac{L_X}{t^i}\circ\cdots\circ\frac{L_X}{t^i}}_{r}(w_i)\=\\
t^{(r+1)(i+1)}L_X\circ\underbrace{\frac{L_X}{t^i}\circ\cdots\circ\frac{L_X}{t^i}}_{r-1}\frac{v_i(t)t^i}{t^i}&\=0\,,
\end{align*}
where the first identity follows from \eqref{eqn:rlk}, the second identity from~\eqref{eqn:intrep} and the definition of~$w_i$, and the last identity again from~\eqref{eqn:rlk} and the assumption on~$v_i(t)$.

The analogous statement in the case~$i=0$, i.e., for solutions of~$L_X^rY=0$, can be proven by using~\eqref{eqn:intrep}.

\smallskip

Proposition~\ref{thm:defde} implies that~$\tfrac{\partial^{r+1}y_\rho}{\partial\rho_i^{r+1}}$ is the sum of a multiple of~$w_i$,  with~$v_i$ in~\eqref{eqn:wi} replaced by $\tfrac{\partial^{r}y_\rho}{\partial\rho_i^{r}}$, and a holomorphic solution of~$L_{X,r(i+1)-1}\circ\cdots\circ L_{X,2i+1}\circ L_{X,i}\circ L_X(u)=~0.$ To prove the corollary we should determine this linear combination. We do it by looking at the coefficients of the local expansion in~$t=0$ of~$w_i(t)$ and~$\tfrac{\partial^{r+1}y_\rho(t)}{\partial\rho_i^{r+1}}$.

A closer look to the recursion formula~\eqref{eqn:rec} for the coefficients of~$y_\rho(t)=\Phi_0(\rho,t)=\sum_{s=0}^\infty{A^{(0)}_s(\rho)t^s}$ reveals that~$A^{(0)}_s(\rho)$ is a polynomial in~$\rho_i$ of degree~$\bigl\lfloor\tfrac{s}{i+1}\bigr\rfloor$ and that the coefficient of~$\rho_i^r$ in~$A^{(0)}_{r(i+1)}(\rho)$ is~$(-1)^{nr}\kappa^{-r}r!^{-2}(i+1)^{-2r}$, where~$\kappa:=\prod_{j=1}^{n-2}{a_j}$, as follows from the explicit expression for~$\hat{P}_{n-2}$ given after~\eqref{eqn:rec}. This implies that the local expansion of~$\tfrac{\partial^{r+1}y_\rho}{\partial\rho_i^{r+1}}$ in~$t=0$ is given by
\begin{equation}
\label{eqn:exp}
\frac{\partial^{r+1}y_\rho(t)}{\partial\rho_i^{r+1}}\=\frac{t^{(r+1)(i+1)}}{(r+1)!^2(i+1)^{2(r+1)}\kappa^{r+1}(-1)^{n(r+1)}}+O(t^{(r+1)(i+1)+1})\,.
\end{equation}
On the other hand, from the definition~\eqref{eqn:wi} with~$v_i(t)=\tfrac{\partial^{r}y_\rho}{\partial\rho_i^{r}}$ is it easy to see that the expansion of~$w_i(t)$ in~$t=0$ is also of the form~\eqref{eqn:exp}.
Since every other solution of~$L_{X,r(i+1)-1}\circ\cdots\circ L_{X,i}\circ L_X(u)\=0$ contains smaller powers of~$t$ in its local expansion, this concludes the proof of the corollary.
\end{proof}

\subsection{Deformation of modular forms}
\label{sec:defmf}
We complete the proof of Theorem~\ref{thm:main}. Recall that the Hauptmodul~$t$ gives the identification~$\Po/\Gamma\simeq\pro^1\smallsetminus\{a_1,\dots,a_{n-2},a_{n-1}=0,a_n=\infty\}$ and that we denote~$P(t):=\prod_{j=1}^{n-1}(t-a_j)$.
\begin{proof}
The result for a modular form~$g\in M_k(\Gamma)$ follows easily from the computation of~$\partial_if$ and~$\partial_it$ since~$\partial_i$ is a derivation and~$g=f^kR(t)$. 
By definition~$f_\rho=y_\rho\bigl(t_\rho(Q)\bigr)$ so then 
\[
\partial_i{f}(q)\=\frac{\partial f_\rho(Q)}{\partial\rho_i}\=\frac{\partial t_\rho(Q)}{\partial\rho_i}\frac{\partial y_\rho(t)}{\partial t}\+\frac{\partial y_\rho(t)}{\partial\rho_i}\circ{t_\rho}\,.
\]
We only need to express the above function as a function of~$Q_\rho=\exp\bigl(\hat{y}_\rho(t)/y_\rho(t)\bigr)$ and then specialize to the Fuchsian value~$\rho=\rho_F$. In order to do that, consider the change of variable formula
\begin{equation}
\label{eqn:changevar}
\frac{1}{Q_\rho}dQ_\rho\=\frac{\kappa}{P(t)y_\rho(t)^2}dt,\qquad\kappa=(-1)^{n}\prod_{j=1}^{n-2}{a_j}.
\end{equation}
This identity follows from well-known properties of the Wronskian~$W(t)$ of the differential equation~\eqref{eqn:de}, namely ~$W(t)=\kappa/P(t)$ and~$y_\rho(t)^2dQ_\rho(t)/dt=W(t)$.

\smallskip

Consider~$\partial y_\rho(t)/\partial\rho_i\circ t_\rho(Q)$. By using Corollary~\ref{thm:integralform} (with $r=0$) and \eqref{eqn:changevar} we get
\[
\frac{\partial y_\rho(t)}{\partial\rho_i}\circ t_\rho(Q)\= f_\rho(Q)\bigintssss_{0}^{Q}\int_0^{Q_1}{h_{\rho,i}(Q_2)\,\frac{dQ_2}{Q_2}\,\frac{dQ_1}{Q_1}}\,,
\]
where~$h_{\rho,i}(Q)=\kappa^{-2}f^4_\rho(Q)t^i_\rho(Q)P(t_\rho(Q))$. 
When we specialize to the Fuchsian value~$\rho_F$ we see that~$h_i(\tau):=h_{\rho_F,i}(q)$ is a weight four cusp form. In fact, it is of weight four since~$f$ is of weight one, it is holomorphic because the pole of~$t^i$ is killed by the zeros of~$f^4$ (see the end of Section~\ref{sec:premod}), and it is zero at every cusp~$c_j$ because~$t(c_j)=a_j,\,j=1,\dots,n,$ up to reordering the indices (a different proof that~$h_i(\tau)$ is a cusp form will be given in Section~\ref{sec:teich}). This proves that~$\partial y_\rho(t)/\partial\rho_i\circ t_\rho(Q)\bigl|_{\rho=\rho_F}=f\widetilde{h}_i'$

\smallskip

The relation~$t_\rho(Q_\rho(t))=t$ implies
\[
\frac{\partial t_\rho(Q)}{\partial\rho_i}\=-Q\frac{\partial t_\rho(Q)}{\partial Q}\cdot\frac{\partial\bigl(\hat{y}_\rho(t)/y_\rho(t)\bigr)}{\partial\rho_i}\circ t_\rho(Q)\,.
\]
Since~$\hat{y}_\rho(t)/y_\rho(t)=\int_{0}^t{\kappa/(y_\rho^2(t_1)P(t_1))\,dt}$, and again using Corollary~\ref{thm:integralform} and formula~\eqref{eqn:changevar}, we finally get
\begin{align*}
\frac{\partial\bigl(\hat{y}_\rho(t)/y_\rho(t)\bigr)}{\partial\rho_i}\circ t_\rho(Q)&\=\int_0^{t_\rho(Q)}\frac{\partial}{\partial\rho_i}\frac{\kappa}{y_\rho^2(t_1)P(t_1)}\,dt_1 \\
&\=-2\int_0^Q\int_0^{Q_1}\int_0^{Q_2}{h_{\rho,i}(Q_3)\,\frac{dQ_3}{Q_3}\,\cdots\,\frac{dQ_1}{Q_1}}\,.
\end{align*}
By specializing to the Fuchsian parameter~$\rho_F$ we find that~$\partial_it=2t'\widetilde{h}_i$ which, combined with~$\partial y_\rho(t)/\partial t\circ t_\rho(Q)|_{\rho=\rho_F}=f'/t'$, concludes the proof.
\end{proof}

The same techniques can be used to show that the higher-order derivatives~$\partial_i^r{g}$ are described by combinations of iterated integrals of Rankin-Cohen brackets of~$g$ and~$h$. As~$r$ grows their modular properties became weaker, but they may still be of some interest (see the third remark in Section~\ref{sec:final}).


\subsection{Teichmüller theory}
\label{sec:teich}
In this section we restrict to the case of modular functions and describe the deformation operators in terms of the deformations of the underlying punctured sphere. This gives an alternative explanation for the appearance of Eichler integrals of weight four cusp forms.
We start by recalling few basic facts about Teichmüller theory.

Let $\Gamma$ be a Fuchsian group such that~$\Po/\Gamma$ is isomorphic to a punctured sphere~$X$ and let~$t\colon\Po/\Gamma\to X$ be a Hauptmodul. Let~$\mathcal{T}(\Gamma)$ denote the Teichm\"uller space of~$\Gamma.$ It is well known that the holomorphic cotangent space at the point~$\Gamma\in\mathcal{T}(\Gamma)$ is the space~$Q(\Gamma)=S_4(\Gamma)$ of quadratic differentials on~$\Po/\Gamma$. There exists a linear isomorphism between~$Q(\Gamma)$ and the space~$D_2(X)$ of rational functions on~$\hat{\C}$ with at most simple poles at the finite punctures of~$X$ and order~$O(|z|^3)$ as~$z\to\infty$ (see Section 2.5 in~\cite{TZ}). This map can be explicitly given in terms of the Hauptmodul~$t$ by
\begin{equation}
\label{eqn:isoratcusp}
R(z)\mapsto q(\tau)\:=R\bigl(t(\tau)\bigr)\cdot t'(\tau)^2\,,\quad z\in\hat{\C},\,\tau\in\Po\,.
\end{equation}

The holomorphic tangent space to $\mathcal{T}(\Gamma)$ at $\Gamma$ is the space $\mathcal{H}(\Gamma)$ of harmonic Beltrami differentials.
The tangent and cotangent spaces are related by the linear map
\begin{equation}
\label{eqn:HQ}
\Lambda^*\colon Q(\Gamma)\to\mathcal{H}(\Gamma),\quad q\mapsto \mathrm{Im}(\tau)^2\bar{q}(\tau)\,\,\quad\tau\in\Po\,.
\end{equation}
Harmonic Beltrami differentials of bounded norm can be used to describe deformations of the punctured sphere~$X$ as follows. Let~$\nu\in \mathcal{H}(\Gamma)$ with~$\|\nu\|_\infty<1$ and denote by~$\mu_\nu\colon\C\to\C$ the measurable function obtained by extending~$\nu$ to~$\C$ by reflection across the real line. Consider the Beltrami differential equation
\begin{equation}
\label{eqn:Belt}
f_{\bar{z}}=\mu_\nu(z)f_{z},\quad z\in\C.
\end{equation}
It is well known that~\eqref{eqn:Belt} has a unique normalized solution~$f^\nu$ that is a homeomorphism of~$\C$ and fixes the points~$0,1,\infty$. The restriction of~$f^\nu$ to~$\Po$ is still a homeomorphism by construction and the conjugate group~$\Gamma^\nu:=f^\nu\Gamma \left(f^\nu\right)^{-1}$ is Fuchsian. It follows that the quotient~$X^\nu:=\Po/\Gamma^\nu$ is a Riemann surface homeomorphic to~$X$ (check Ahlfors's book~\cite{ahlfors} for more details and proofs of these statements). 

\smallskip

As the above paragraph shows, one may construct deformations of~$X$ starting with a rational function in~$D_2(X)$ by composing the maps in~\eqref{eqn:isoratcusp} and~\eqref{eqn:HQ} (but this map in general do not give Beltrami differentials with the required norm). The coefficient of the accessory parameter~$\rho_i$ in the differential equation~\eqref{eqn:de} in normal form is~$R_i(t)=t^i/P(t)\in D_2(X)$. The quadratic differential associated to~$R_i$ via the map~\eqref{eqn:isoratcusp} is precisely the cusp form~$h_i(\tau)$ appearing in Theorem~\eqref{thm:main}, as follows from the identity~$t'=P(t)f^2\kappa^{-2}$. In this way, using~\eqref{eqn:HQ}, we can associate to every accessory parameter~$\rho_i,\,i=0,\dots,n-3$ a harmonic Beltrami differential~$\nu_i:=\Lambda^*(h_i(\tau))$. Moreover,~$\nu_0,\dots,\nu_{n-4}$ form a basis of~$\mathcal{H}(\Gamma)$.

Now fix~$0\le i\le n-4$ and let $0\neq\vare\in\C$ be such that~$\|\vare\nu_i\|_\infty<1$. Denote by~$f^{\nu_i}$ the normalized homeomorphic solution of the Beltrami equation~\eqref{eqn:Belt} with~$\mu_\nu=\mu_{\vare\nu_i}$. 
As above, one obtains a Fuchsian group~$\Gamma_i:=f^{\vare\nu_i}\Gamma(f^{\vare\nu_i})^{-1}$ and an~$n$-punctured sphere~$X_i:=\Po/\Gamma_i$. The situation, for every~$i=0,\dots,n-4$, is summarized in the following diagram
\begin{equation}
\label{eqn:diag2}
\begin{CD}
\Po @>f^{\vare\nu_i}>> \Po \\
@V{t}VV       @VV{t^{\vare\nu_i}}V \\
X    @>>F^{\vare\nu_i}>   X_i\\
\end{CD}
\end{equation}
where~$t^{\vare\nu_i}\colon\Po\to X_i$ is a Hauptmodul for~$\Gamma_i$ normalized by~$t^{\vare\nu_i}(\infty)=t(\infty),\,t^{\vare\nu_i}(0)=t(0)$ as follows from the normalization of~$f^{\vare\nu_i}$.

It is known (\cite{ahlfors},\cite{TZ}) that~$F^{\vare\nu_i}$ is a quasiconformal map of Riemann sufaces and is holomorphic in~$\vare,$ while both~$f^{\vare\nu_i}$ and~$t^{\vare\nu_i}$ are real-analytic non holmorphic functions in~$\vare$. In particular, it makes sense to consider the derivatives of the above functions with respect to~$\vare$ and~$\bar{\vare}.$ From a well-known formula of Ahlfors (see formulae 3.9 and 4.3 in~\cite{TZ})  it follows that, for~$0\le i\le n-4$,
\begin{equation}
\label{eqn:defteich}
\frac{\partial{t^{\vare\nu_i}}}{\partial\bar{\vare}}\Bigr|_{\vare=0}\=\frac{1}{2}\widetilde{h}_it\,.
\end{equation}
The right hand-side of~\eqref{eqn:defteich} reminds the statement Theorem~\ref{thm:main} in the case of a modular function (weight k=0). The reason is the following.
\begin{theorem}
\label{thm:def}
Let $X$ be an $n$-punctured sphere, and let $t\colon\Po/\Gamma\to X$ be a Hauptmodul.  
For~$i=0,\dots,n-4$, let $\nu_i=\Lambda^*(h_i)\in\mathcal{H}(\Gamma)$ and let~$\partial_i$ be the~$i$-th deformation operator on~${M}_*(\Gamma)$. Then
\[
\partial_i t\=4\,\frac{\partial{t^{\vare\nu_i}}}{\partial\bar{\vare}}\Bigr|_{\vare=0}.
\]
\end{theorem}
\begin{proof}
Let $X=\pro^1\smallsetminus\{a_1,\dots,a_n=\infty\}$. Fix~$0\le i\le n-4$, let~$0\neq\vare\in\C$ be such that~$\|\vare\nu_i\|_\infty<1$ and consider the $n$-punctured sphere~$X_i=\pro^1\smallsetminus\{a^{\vare\nu_j}_1,\dots,a^{\vare\nu_j}_n=\infty\}$ where $a_j^{\vare\nu_i}:=F^{\vare\nu_i}(a_j),\,j=1,\dots,n$ (see~\eqref{eqn:diag2}).

To the Fuchsian uniformization of $X_i$ is associated a differential equation~\eqref{eqn:de} with singular points~$a_j^{\vare\nu_i}$ and accessory parameters~$\rho_0^{\vare\nu_i},\dots,\rho_{n-4}^{\vare\nu_i}.$ These accessory parameters are continuously differentiable in $\vare$ since they are coefficients of the $q$-expansion of $t^{\vare\nu_i}$ and this function is real-analytic in $\vare$. 
The theorem is a consequence of the identity 4.6 of~\cite{TZ}, namely\begin{footnote}{The identity we refer to in~\cite{TZ} is stated for accessory parameters~$c_1,\dots,c_{n}$ of a differential equation projectively equivalent to~\eqref{eqn:de}. A straightforward computation shows that those are related to the accessory parameters~$\rho_0,\dots,\rho_{n-4}$ by~$c_i=\Res_{t=\alpha_i}\bigl(4\sum_{i=0}^{n-3}{\rho_it^i}-P''(t)\bigl)/2P(t)$. This leads to the identity in~\eqref{eqn:claim}.}
\end{footnote}
\begin{equation}
\label{eqn:claim}
\frac{\partial{\rho_j^{\vare\nu_i}}}{\partial\bar{\vare}}\Bigr|_{\vare=0} =
\begin{cases} 
\frac{1}{4}& \quad i=j\,,\\
0 &\quad i\neq j\,.
\end{cases}
\end{equation}

The reason is the following. Let $t^{\vare\nu_i}(\tau)=\sum_{s=1}^{\infty}{t_s^{\vare\nu_i}}q^s$ be the $q$-expansion at $\infty$ of the normalized Hauptmodul $t^{\vare\nu_i}$ of~$\Gamma_i$. If the~$q$-expansion of the Hauptmodul~$t$ of~$\Gamma$ is~$t(\tau)=\sum_{s=1}^\infty{t_s(\rho,a)q^s}\,,a=(a_1,\dots,a_{n-1})$, then the Fourier coefficients of~$t^{\vare\nu_j}$ are of the form
\[
t_s^{\vare\nu_j}=t_s(\rho^{\vare\nu_j},a^{\vare\nu_j})\,,\quad a^{\vare\nu_j}=(a_1^{\vare\nu_j},\dots,a_{n-1}^{\vare\nu_j})\,.
\]
Now, $a_j^{\vare\nu_i}$ is holomorphic in $\vare;$ this follows from the definition $a_j^{\vare\nu_i}=~F^{\vare\nu_i}(a_j)$ and the fact that $F^{\vare\nu_i}$ is holomorphic in $\vare.$ This implies that the derivative of $\alpha_j^{\vare\nu_i}$ with respect to $\bar{\vare}$ is zero and then
\begin{equation}
\label{eqn:Ft}
\frac{\partial{t^{\vare\nu_i}}}{\partial\bar{\vare}}{\Bigr|_{\vare=0}}\= \sum_{s=1}^\infty{\left(\sum_{k=0}^{n-4}{\frac{\partial t_s(\rho^{\vare\nu_i},a^{\vare\nu_i})}{\partial{\rho_k^{\vare\nu_i}}}\Bigr|_{\vare=0}\frac{\partial\rho_k^{\vare\nu_i}}{\partial\bar{\vare}}\Bigr|_{\vare=0}}\right)q^s}.
\end{equation}
On the other hand, the action of $\partial_i$ on $t$ is, by definition, 
\begin{equation}
\label{eqn:Ftj}
\partial_it\=\sum_{s=1}^\infty{\frac{\partial t_s(\rho,a)}{\partial{\rho_i}}q^s}\,.
\end{equation}
By comparing the~$q$-expansions~\eqref{eqn:Ft} and~\eqref{eqn:Ftj} and using the identity
\[
\frac{\partial t_s(\rho^{\vare\nu_i},a^{\vare\nu_i})}{\partial{\rho_k^{\vare\nu_i}}}\Bigr|_{\vare=0}\=\frac{\partial t_s(\rho,a)}{\partial{\rho_k}}\Bigr|_{\rho=\rho_F}\,,\quad s\ge1,\,k=0,\dots,n-4\,,
\]
together with~\eqref{eqn:claim}, the statement of the theorem follows.
\end{proof}
The above proposition together with~\eqref{eqn:defteich} gives another proof of Theorem~\ref{thm:main} in the case of (meromorphic) modular forms of weight zero.

\subsection{Vector-valued modular forms}
\label{sec:vvmf}
In this final section we reformulate the results of~\ref{sec:defde} and \ref{sec:defmf} in terms of vector-valued modular forms. From this perspective, we shall consider two situations: the action of~$\partial/\partial\rho_i$ on the space of solutions of~$L_X$ and the action of~$\partial_i$ on~$M_*(\Gamma)$ induce different maps between spaces of modular forms and VVMFs attached to certain extensions. 

In Section~\ref{sec:ext} we showed that extensions of symmetric tensor representations can be described in terms of quasimodular forms. Let~$h\in S_4(\Gamma)$ and let~$p_h(\gamma,\tau):=r_{h,2}(\gamma)\tau^2+r_{h,1}(\gamma)\tau+r_{h,0}(\gamma)$ be its period polynomial. From~\eqref{eqn:ext} it follows that~$h$ induces extensions~$[v^h_{0,2}]\in\Ext^1_\Gamma(v_0,v_2)$ and $[v^h_{1,1}]\in\Ext^1_\Gamma(v_1,v_1)$. We can describe explicitly a representative of each class in terms of~$p_h(\gamma,\tau)$ as follows
\begin{equation}
\label{eqn:v02}
\gamma=\begin{pmatrix}
a & b \\
c & d
\end{pmatrix}\;\;\mapsto\;\;v^h_{0,2}(\gamma)=
\left(\begin{array}{c|ccc}
1 & r_{h,2}(\gamma) & r_{h,1}(\gamma) & r_{h,0}(\gamma) \\
\hline
0 & a^2 & 2ab & b^2 \\
0 & ac & ad+bc & bd \\
0 & c^2 & 2cd & d^2
\end{array}\right),
\end{equation}
 \[
\gamma=\begin{pmatrix}
a & b \\
c & d
\end{pmatrix}\;\;\mapsto\;\;v^h_{1,1}(\gamma)=
\left(\begin{array}{c|c}
\begin{matrix}
a & b \\
c & d
\end{matrix} & \gamma\cdot B_h(\gamma)\\
\hline
\begin{matrix}
0 & 0 \\
0 & 0
\end{matrix} & 
\begin{matrix}
a & b \\
c & d
\end{matrix}
\end{array}\right),\quad B_h(\gamma):=\begin{pmatrix}
r_{h,1}(\gamma) & -2r_{h,0}(\gamma) \\
2r_{h,2}(\gamma) & -r_{h,1}(\gamma)
\end{pmatrix}.
\]

The action of~$\partial/\partial\rho_i$ on the space of solutions of~$L_X$ is related to extensions~$v_{1,1}^{h_i}$ of two-dimensional symmetric tensor representations of~$\Gamma$.
\begin{proposition}
\label{thm:vv0}
For every~$i=0,\dots,n-4,\,r\ge0,$ the derivative~$(\partial^r y_\rho/\partial\rho_i^r)$ lifts to a component of a vector-valued modular form with respect to a~$r$-iterated extension of symmetric tensor representations of dimension~$2$. 
When~$r=1$ the derivation~$\partial/\partial\rho_i$ on~$y_\rho(t)$ induces a map
\[
M_k(\Gamma)\to M_k(\Gamma,v_{1,1}^{h_i}),\quad g\mapsto\begin{pmatrix}
\tau{g}\widetilde{h}'_i-g\widetilde{h}_i \\
g\widetilde{h}_i'\\
\tau{g}\\
g
\end{pmatrix}\,,
\]
for every~$k\ge 0$, where~$h_i\in S_4(\Gamma)$ is as in Theorem~\ref{thm:main}.
\end{proposition}
\begin{proof}
The result follows from Proposition~\ref{thm:defde}. Recall that~$\partial^ry_\rho(t)/\partial\rho_i^r$ is a solution of a Fuchsian equation obtained as the composition of~$(r+1)$ second-order Fuchsian operators. 
The vector of solutions of this equation composed with~$t(\tau)$ gives a vector of holomorphic functions on~$\Po$. Its transformation property follows from the fact that the monodromy of the differential equation is an iterated extension of symmetric tensor representations. 

In the case~$r=1$, the proof of Theorem~\ref{thm:main} shows that~$\partial y_\rho(t)/\partial\rho_i$ lifts to~$f\widetilde{h}_i'$ and a similar computation shows that~$\partial\hat{y}_\rho(t)/\partial\rho_i$ lifts to $\tau{f}\widetilde{h}'_i-f\widetilde{h}_i$. The modular transformation properties of these functions show that the extension has to be~$v_{1,1}^{h_i}$. 
\end{proof}

The deformation operator~$\partial_i$ is related to extensions in~$\Ext^1_\Gamma(v_0,v_2)$. 
\begin{proposition}
\label{thm:vv}
For every~$i=0,\dots,n-4$ the~$i$-th deformation operator~$\partial_i$ induces a map
\[
\overset{\rightarrow}{\partial}_i\colon M_k(\Gamma)=M_k(\Gamma,v_0)\to M_k(\Gamma,v^{h_i}_{0,2})\quad g\mapsto\overset{\rightarrow}{\partial}_i{g}\:=\begin{pmatrix}
\partial_i{g} \\
\tau^2{g'}+2\tau{g}\\
\tau g'+g\\
g'
\end{pmatrix}\,,
\] for every~$k\ge 0$, where~$h_i\in S_4(\Gamma)$ is as in Theorem~\ref{thm:main}.
\end{proposition}
\begin{proof}
To check that~$\overset{\rightarrow}{\partial_i}{g}$ is a VVMF with respect to~$v_{0,2}^{h_i}$ we split the vector into the lower part~$(g',\tau g'+g,\tau^2g'+\tau g)^{t}$ and the upper part~$(\partial_i{g})$ and check that they transform accordingly under the action of~$\Gamma$.

The vector~$(g',\tau g'+g,\tau^2g'+\tau g)^{t}$ is a weight~$k$ VVMF for the symmetric tensor representation~$v_2$ associated to the quasimodular form~$g'$, as Choie-Lee's paper~\cite{CL} or a simple check shows. 
On the other hand, from Theorem~\ref{thm:main} it follows that
\[
\partial_i{g}\bigl|_k\gamma\=[g,\widetilde{h}_i]\bigl|_k\gamma\=[g,\widetilde{h}_i]\+[g,p_{h_i}(\gamma,\tau)]\,,
\]
and an easy computation shows~$[g,p_{h_i}(\gamma,\tau)]=r_{h_i,2}(\tau^2g+2\tau{g})+r_{h_i,1}(\tau{g'}+g)+r_{h_i,0}g'$. By comparing these transformations with the explicit description of~$v_{2,0}^{h_i}$ in~\eqref{eqn:v02} the statement follows.
\end{proof}

\subsection{Final remarks}
\label{sec:final}
\begin{enumerate}
\item It should be possible, by considering differential equations of higher order satisfied by~$g'$ and~$\partial_i{g}$, to prove Proposition~\ref{thm:vv} by a monodromy argument like Proposition~\ref{thm:vv0}. The reason why the same argument does not work is related to the appearance of~$g'$ in~$\overset{\rightarrow}{\partial}_i{g}:$ the non-trivial depth of this quasimodular form does not permit to reduce to second order differential equations as happens for modular forms.
\item A related problem is to extend the map~$\overset{\rightarrow}{\partial}_i$ to a map on~$M_*(\Gamma,v_l)$ for every~$l>0$. It follows from~\eqref{eqn:cl} that this is equivalent to define the deformation operators on quasimodular forms.  As quasimodular forms does not fit into the classical picture of uniformizing differential equations, they are not a priori related with the accessory parameters. Nevertheless, by writing quasimodular forms as weight zero VVMFs, one can argue as in Chapter 5 of~\cite{123} and find differential equations solved by quasimodular forms. If these differential equations are special members of a family depending on some parameters, then one can define deformations with respect to those parameters. This may give a reasonable way to extend the maps~$\overset{\rightarrow}{\partial}_i$.
\item By using Corollary~\ref{thm:integralform} it is possible to compute higher deformations~$\partial_i^r{g}$ as well as mixed derivatives in terms of iterated integrals of modular forms. More generally, given a genus zero group~$\Gamma$ with~$n$-cusps it may be interesting to consider the expansion of a modular form~$f\in M_*(\Gamma)$ around the Fuchsian value of the accessory parameter
\[
\hat{f}(\tau)\=\sum_{J\in\Z^{n-4}}^{\infty}{\partial^J{f(\tau)}(\rho-\rho_F)^J},\quad \partial^J:=\partial_1^{j_1}\cdots\partial_{n-4}^{n-4}\,,
\]
and to investigate its modular properties with respect to the usual action of~$\Gamma$.
\end{enumerate}

\section*{Acknowledgments}
The paper was written while I was a graduate student at SISSA (International School for Advanced Studies) in Trieste and a visiting student of the IMPRS graduate school in the Max Planck Institute for Mathematics in Bonn. I want to thank both institutions for the excellent working conditions. I want to thank my advisers Don Zagier and Fernando Rodriguez Villegas for their suggestions and support, and Vasily Golyshev for very useful conversations.


\begin{thebibliography}{20}
\bibitem{ahlfors}{L.\,Ahlfors}, `Lectures on quasiconformal mappings', with additional chapters by C.J.\,Earle\, and I.\,Kra, M.\,Shishikura, and J.H.\,Hubbard, {\em University Lecture Series}, American Mathematical Society 38 (2006).
%
\bibitem{beukers}{F.\,Beukers}, `Irrationality proofs using modular forms', Journées Arithmétiques de Besancon, {\em Astérisque} 147-148. (1987), 271-283.
%
\bibitem{bogo}{G.\,Bogo}, `Accessory parameters for four-punctured spheres', {\em Arxiv:2004.02971}.
%
\bibitem{123}{J.\,Bruinier, G.\,Harder, G.\,van der Geer, D.\,Zagier}, `The 1-2-3 of Modular Forms: Lectures at a Summer School in Nordfjordeid, Norway' (ed. K. Ranestad), {\em Universitext}, Springer-Verlag, Berlin-Heidelberg-New York (2008).
%
\bibitem{CL}{Y.\,Choie, MH.\,Lee}, `Symmetric tensor representations, quasimodular forms and weak Jacobi forms', {\em Advances in Mathematics}, 287 (2016) 567-599.
%
\bibitem{KS}{M.\,Kuga, G.\,Shimura}, `On vector differential forms attached to automorphic forms', {\em J. Math. Soc. Japan} 12 3 (1960) 258-270.
%
\bibitem{MR}{M.\,Mertens, M.\,Raum}, `Modular forms of virtually real-arithmetic type I: Mixed mock modular forms yield vector-valued modular forms', {\em Mathematical Research Letters} 28 (2021), 511-561. 
%
\bibitem{sgervais}{H.\,P.\,de Saint-Gervais}, `Uniformization of Riemann Surfaces: revisiting a hundred-year-old theorem', \emph{European Mathematical Society} (2016).
%
\bibitem{TZ}{L.\,Takhtajan, P.\,Zograf}, `On the Liouville equation, accessory parameters and the geometry of the Teichmüller space for the Riemann surfaces of genus 0', {\em Mathematics of the USSR-Sbornik} 60 (1988) 143-161.
\end{thebibliography}
\end{document}